\documentclass[12pt]{article}

\usepackage{amsmath,amssymb,amsfonts,bm,slashed,graphics}
\usepackage[latin1]{inputenc}
\usepackage[english]{babel}
\usepackage{epsfig}
\usepackage{hyperref}
\usepackage{amsthm}
\usepackage{subfigure}
\usepackage{dsfont}
\usepackage{hyperref}
\usepackage{amsthm}
\usepackage{subfigure}
\usepackage{color}
\usepackage{multirow}
\usepackage{psfrag}
\usepackage{graphicx}
\usepackage{extarrows}
\usepackage{color}
\usepackage{multirow}
\usepackage{tikz}

\DeclareMathOperator{\Tr}{Tr}

\newcommand{\cT}{{\cal{T}}}

\newtheorem{lemma}{Lemma}

\newtheorem{definition}{Definition}
\newtheorem{theorem}{Theorem}
\newtheorem{proposition}{Proposition}

\newcommand{\bea}{\begin{eqnarray}}
\newcommand{\eea}{\end{eqnarray}}
\newcommand{\bee}{\begin{equation}}
\newcommand{\ee}{\end{equation}}

\def\C{{\mathbf C}}

\newcommand{{\J}}{{\mathbf J}}
\newcommand{{\bbR}}{{\mathbb R}}
\newcommand{{\bbJ}}{{\mathbb J}}
\newcommand{{\bbL}}{{\mathbb L}}
\newcommand{{\bbB}}{{\mathbb B}}
\newcommand{{\bbS}}{{\mathbb S}}
\newcommand{{\bbZ}}{{\mathbb Z}}
\newcommand{{\bbT}}{{\mathbb T}}

\newcommand{{\bbw}}{{\mathbb w}}

\newcommand{{\frw}}{{\mathfrak  w}}
\newcommand{{\frW}}{{\mathfrak  W}}

\newcommand{{\mcT}}{{\mathcal T}}

\newcommand{\cR}{{\cal R}}

\newcommand{\cH}{{\cal H}}

\newcommand{\cC}{{\cal C}}

\newcommand{\cF}{{\cal F}}

\newcommand{\bbone}{{\bf 1}}

\newcommand{\beann}{\begin{eqnarray*}}
\newcommand{\eeann}{\end{eqnarray*}}  

\newcommand{\gF}{{\mathfrak{F}}}

\newcommand{\prf}{{\noindent {\bf Proof}\quad }}
\begin{document}
\title{Loop Vertex Representation for Cumulants,\\
Part II: Weingarten Calculus}
\author{V. Rivasseau\\
Universit\'e Paris-Saclay, CNRS/IN2P3\\ IJCLab, 91405 Orsay, France\\
Email: vincent.rivasseau@gmail.com}
\date{} 

\maketitle

\begin{abstract}  
 In this paper we construct scalar cumulants for stable random matrix models with single trace interactions of arbitrarily 
 high even order by Weingarten calculus.  We obtain explicit and convergent expansions for these scalar cumulants
in the limit $N\to \infty$.
%We show that any cumulant is an analytic function inside a cardioid domain  in the complex plane and we  
%prove their Borel-LeRoy summability at the origin of the  coupling constant. 
%Our proof is uniform in the external variables. 
\end{abstract}

\noindent\textbf{keywords}
Random Matrix; Cumulants; Constructive Field Theory


\section{Introduction}

We are interested in the loop vertex representation  (LVR) \cite{rivasseau2018loop,KRS,KRS1}, which is a improvement of the loop vertex expansion  (LVE) \cite{rivasseau2007constructive}. This LVE was introduced itself as a tool for constructive field theory to deal with random matrix fields. 

\medskip
First, constructive field theory 
\cite{simon2015p,glimm2012quantum,rivasseau2014perturbative} focused at the interaction and the parameter 
$\lambda$ which measures the strength of this interaction. Naturally, we set about solving the problem. A common main feature of the LVR and LVE is that it is written in terms of trees which are exponentially bounded. It means that the outcome of the LVR-LVE is convergent and is \emph{the Borel-LeRoy sum in $\lambda$},
whereas the  usual perturbative quantum field theory \emph{diverges} at the point $\lambda=0$. 
The essential components of LVE are the Hubbard-Stratonovich intermediate field
representation \cite{Hub,Str}, the replica method \cite{MPV} and the BKAR formula \cite{BK,AR1}.
The added ingredients of the LVR are combinatorial, based on the selective Gaussian integration \cite{rivasseau2018loop}, and the Fuss-Catalan numbers and their generating function \cite{FussCatalan}.
We think that the LVR has \emph{more power} than the LVE, since the LVR can treat more models, 
with higher polynomial interactions\footnote{For an early application to the generating function
of connected Schwinger functions - which in this paper are denoted {\it cumulants} - see \cite{magnen2008constructive}; 
for the actual mechanism of replacing Feynman graphs, which are not exponentially bounded, by trees, see \cite{RiZh1,RiTa};
for a review of the LVE, we suggest consulting \cite{GRS}.}.

 \medskip
This article is devoted to Weingarten calculus \cite{Weingarten,Col,ColSni} and their associated \emph{combinatorial maps}, and is a sequel to \cite{GuKra,Riv1}. Application of random matrices to 2d quantum gravity \cite{matrix} relies on their combinatorial maps, which depend on (at least) two parameters: a coupling constant $\lambda$ and the size of the matrix, $N$. A \emph{formal} expansion in the parameter $\lambda$ yields generating functions for maps of arbitrary genus. The coupling constant $\lambda$ roughly measures the size of the map while the parameter $1/N$ turns out to measure the genus of the map \cite{thooft}. 

 \medskip
Our real purpose is however to study the distributions \emph{asymptotically}, when the size of the matrices or the tensors $N$ goes to infinity. For classical random matrices, 
 the \emph{combinatorial maps of genus $0$} constitute the order $0$ in $N$.
 For the random tensors, Gurau and coworkers
 discovered that the order $0$ in $N$ is a subset of combinatorial maps,  
 and they  called it \emph{melonic} 
 \cite{1Nexpansion1, 1Nexpansion2, 1Nexpansion3, critical, Gurau-universality, uncoloring, bonzom-SD, Gurau-book, enhanced-1, Walsh-maps, multicritical, bonzom-review, Lionni-thesis, 
 Gurau:2012ix, Gurau:2013pca}. 

\medskip
The authors of \cite{KRS}  join the LVR to Cauchy holomorphic matrix calculus
and have been applied to the simplest complex matrix model with stable monomial interaction. In \cite{KRS1} the  same authors have extended it to the case of \emph{Hermitian} or \emph{real symmetric} matrices,
in a manner both \emph{simpler and more powerful}. The basic formalism is still the LVR, but
while \cite{rivasseau2018loop,KRS} used contour integral parameters attached to every \emph{vertex} 
of the loop representation, \cite{KRS1} introduces more contour integrals, one for each \emph{loop vertex corner}.
This results in simpler bounds for the norm of the corner operators. 

 \medskip
\begin{figure}[!htb]\centering
\includegraphics[width=0.3\linewidth]{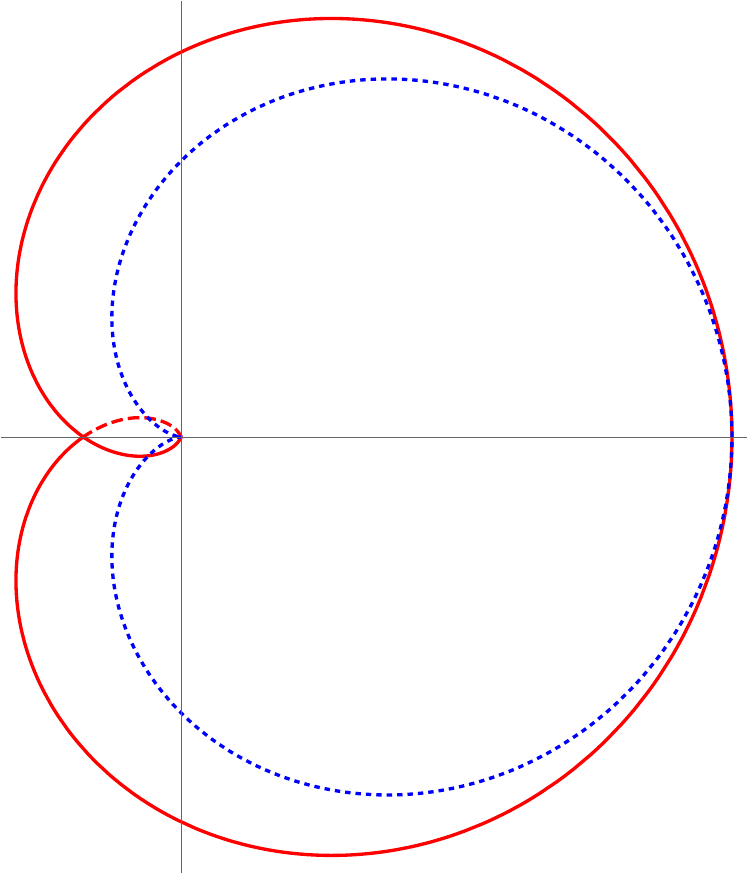} 
\caption{In blue the cardioid domain considered in \cite{rivasseau2022cumulants}, in red the cardioid domain considered
in \cite{BGKL}.}
\label{fig:W-domain}
\end{figure}

\medskip
But we should remember that the LVE is older and their authors have more time 
to fine-tune their models. They construct their models with the coupling constant in a cardioid-shaped domain
(see Figure \ref{fig:W-domain}) which has opening angle arbitrarily close to $2 \pi$  \cite{rivasseau2022cumulants} 
or even exceeding $2 \pi$ \cite{BGKL}. In this case the LVE is capable to compute  
 some typically non-perturbative effects like instantons by resuming perturbative field theory. In \cite{sazonov2023}, Sazonov combined the LVE with ideas of the variational perturbation theory.

\medskip\noindent
{\bf Acknowledgement} 
We acknowledge the support of the IJClab and the CEA-List.

\section{The model} 

In this paper, $\cH$ is the Hilbert space $\cH=C^N$,  
$\Tr $ always means \emph{the trace on $\C^{N}$}, $\Tr_\otimes$ always  means \emph{the trace on $\C^{N\times N}$}, and $\bbone_\otimes$ always  means  the $N^2 \times N^2$ matrix whose all eigenvalues are $1$. 

Consider a complex square matrix model with stable interaction of order $2p$, where $p\ge 2$ is an integer 
which is fixed through all this paper.
We assume the reader is reasonably familiar with the notations of \cite{KRS,KRS1,GuKra}
and with Appendix B of the book \cite{guruau2017random}.
Let us recall some basics of our LVR  in the scalar and $d=0$ case \cite{rivasseau2018loop}. 
One of the key elements of the
LVR construction is the Fuss-Catalan numbers of order $p$, which we denote by $C_{n}^{(p)}$, and their generating function $T_p$ \cite{FussCatalan}.
This generating function $T_p$ is defined by
\begin{equation}
T_p(z) = \sum_{n=0}^\infty C_{n}^{(p)} z^n . \label{gencat1}
\end{equation}
It is analytic at the origin and obeys the algebraic equation 
\begin{equation}\label{gencatalan}
zT_p^p(z) -T_p(z) +1 =0 .
\end{equation}

\medskip
In the case $p=3$ the LVR is somewhat simplified; the Fuss-Catalan equation is 
\bee
zT_3^3(z) -T_3(z) +1 =0 , \label{cardanomain}
\ee
which is soluble by radicals. We give in \cite{rivasseau2018loop}, section VI.2,  the details derived from Cardano's solution.

We shall only present our main result for \emph{complex square matrices} in a perturbation $(M M^\dagger )^p $.
In a simplification with respect to \cite{KRS}, we consider only square matrices.
The generalisation to other cases, for instance rectangular complex matrices, or Hermitian matrices, or real symmetric matrices, 
is not too difficult for someone who is familiar of \cite{KRS,KRS1}.
\footnote{For practical applications such as data analysis, the case $p=3$ seems  to be the main one and it is
interesting  to treat the case of \emph{real symmetric} matrices and  \emph{rectangular matrices}.} 

Next we shall define a mathematical expression for the cumulants.     

\begin{definition}
\label{cumulantsdef}
The cumulant of order $2{\mathcal K}$ is:
\bea
\hskip-.6cm\label{cum01}
\mathfrak{K}^{{\mathcal K}}(\lambda,N )\!&:=&\!\Big[
\frac{\partial^{2}}{J^{\ast}_{a_{1}b_{1}} J_{c_{1}d_{1}}}\cdots
\frac{\partial^{2}}{J^{\ast}_{a_{{\mathcal K}}b_{{\mathcal K}}} J_{c_{{\mathcal K}}d_{{\mathcal K}}}}
\log {\cal Z}(\lambda, N,J)\Big]_{J=0}, 
\eea
where
\bea \label{diffnor}
\log {\cal Z} (\lambda, N,J)= \frac{1}{N^2} \log Z(\lambda, N,J) 
\eea
and $J_{ab}^{\ast}$ is the complex conjugate of $J_{ab}$, so that $(J^{\dagger})_{ab}=J^{\ast}_{ba}$.
\end{definition}

Next we introduce Weingarten functions who were defined in \cite{Col},  \cite{ColSni} and also in \cite{GuKra}.
As the authors of \cite{GuKra} remark, scalar cumulants
arise when integrating over unitary matrices $\text{U}(N)$ with the invariant normalized Haar measure. 
Denoting $U_{ab}^*$ the complex conjugate of $U_{ab}$ we have \cite{Col,ColSni}:
\bea
&&\int dU \; U_{a_{1}b_{1}}\dots U_{a_{k}b_{k}}
U^{*}_{c_{1}d_{1}}\dots U^{*}_{c_{l}d_{l}}=
\\&& \nonumber\hskip2cm
 \delta_{kl}\sum_{\sigma,\tau\in \mathfrak{S}_{k}}
\delta_{a_{\tau(1)c_{1}}}\dots \delta_{a_{\tau(k)}c_{k}}
\delta_{b_{\sigma(1)}d_{1}}\dots \delta_{b_{\sigma(k)}d_{k}}
{\cal{W}}\hskip-.06cm_{g}(\tau\sigma^{-1},N) \; .
\label{Weingartenrelation}
\eea
The functions ${\cal{W}}\hskip-.06cm_{g}(\zeta=\tau\sigma^{-1},N) $ only depend on the cycle structure of $ \zeta$.
Here is a few examples of Weingarten functions:
\begin{align*}
{\cal{W}}\hskip-.06cm_{g}\big((1),N\big)&=\frac{1}{N}& {\cal{W}}\hskip-.06cm_{g}\big((1,1,1),N\big)&=\frac{N^{2}-2}{N(N^{2}-1)(N^{2}-4)}\\
{\cal{W}}\hskip-.06cm_{g}\big((1,1),N\big)&=\frac{-1}{N^{2}-1}& {\cal{W}}\hskip-.06cm_{g}\big((1,2),N\big)&=\frac{-1}{(N^{2}-1)(N^{2}-4)}\\
{\cal{W}}\hskip-.06cm_{g}\big((2),N\big)&=\frac{-1}{N(N^{2}-1)}& {\cal{W}}\hskip-.06cm_{g}\big((3),N\big)&=
\frac{2}{N(N^{2}-1)(N^{2}-4)} \; .
\end{align*}

For any permutation of $k$ elements $\zeta\in{\mathfrak S}_{{k}}$, let us write $C(\zeta)$ the integer partition of 
associated to the cycle decomposition 
of $\zeta$ and $|C(\zeta)|$ the number of cycles it contains. Let us also denote by $\Pi_{{k}}$ the set of integer partitions of $k$ 
(recall that a partition $\pi\in\Pi_{{k}}$  is an increasing sequence of $|\pi|$ integers $0<{k}_{1}\leq \cdots\leq {k}_{|\pi|}$ such 
that ${k}_{1}+ \cdots + {k}_{|\pi|}={k}$).
To any integer partition of ${k}$ we associate a \emph{trace invariant}:
\begin{equation}
\Tr_{\pi}({\mathfrak X})=\Tr({\mathfrak X}^{{k}_{1}})\cdots\Tr({\mathfrak X}^{{k}_{p}}) \; .
\end{equation}

Let us chose a permutation $\zeta\in{\mathfrak S}_{{k}}$ whose 
cycle decomposition reproduces the contribution 
of the \emph{broken faces} to the amplitude of a graph. Specifically, if there are $b=b(G)$ broken faces with ${k}_{1},\dots,{k}_{b}$ cilia,
we choose $\zeta$ to have a cycle decomposition of the form:
\begin{equation}
\zeta=(i_{1}^{1}\dots i_{{k}_{1}}^{1}) \cdots (i_{1}^{b}\dots i_{{k}_{b}}^{b}) \; .
\end{equation}
This permutation defines a labeling of the cilia in such a way that the product of traces over the broken faces can be expressed as:
\bea\label{CilTra}&&
\prod_{1\leq m\leq b}\Tr\Big[JJ^\dagger
\mathop{\prod}\limits_{1\leq r\leq {k}_{m}}^{\longrightarrow}{\mathfrak X}^{i^{m}_{r}}\Big]=
\sum_{1\leq p_{1}, q_1 \dots \leq N}  (JJ^{\dagger})_{p_{l}q_{l}}   \prod_{1\leq l\leq {k}}{\mathfrak X}^{l}_{q_{l}p_{\zeta(l)}} 
\eea
where ${\mathfrak X}^{l}$ is the product of the resolvents  
\bea
\big[\bbone_\otimes   +\Sigma (\lambda,M)\big]^{-1} &=& \big[\bbone_\otimes   +
\lambda \sum_{{\mathfrak k} = 0}^{p-1} A^{\mathfrak k}(MM^\dagger) \otimes A^{p-1-{\mathfrak k}}(M^\dagger M) \big]^{-1}
\label{Xnew}
\eea
located on the corners separating the cilia labeled $l$ and $\zeta(l)$. 

Similarly, for the $f(G)-b(G)$ unbroken faces we denote by ${\mathfrak X}^{m}$ the product of the resolvents 
\bea
\big[\bbone_\otimes   +\Sigma (\lambda,M)\big]^{-1} 
=\big[\bbone_\otimes   +
\lambda \sum_{{\mathfrak k} = 0}^{p-1} A^{\mathfrak k}(MM^\dagger) \otimes A^{p-1-{\mathfrak k}}(M^\dagger M) \big]^{-1}
= {\mathfrak X}^{m}
\label{Ynew}
\eea around the unbroken face labeled $m$.

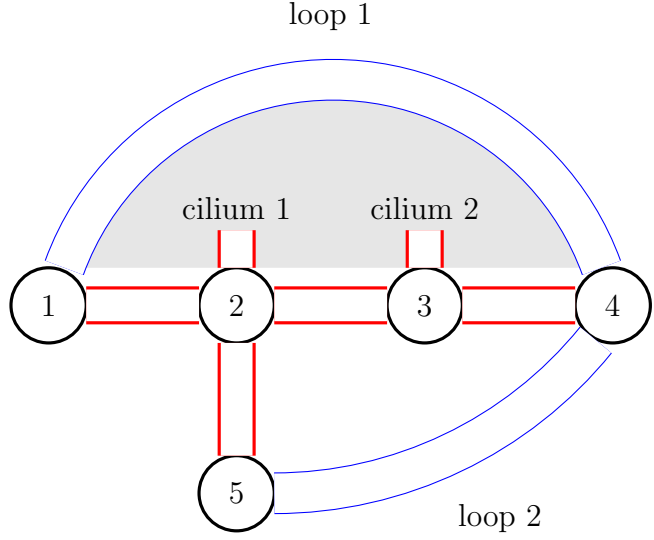
\begin{figure}[!htb]
\centering
\begin{tikzpicture}[>=latex,scale=0.5]

\filldraw[draw=black,fill=gray!20]
plot (14.7,1) arc (20: 160: 7.6)
plot (14,0) -- (11,0) 
-- cycle;
\draw[color=black, very thick](0,0) circle (1);   
\draw[color=black, very thick](5,0) circle (1); 
\draw[color=black, very thick](10,0) circle (1); 
\draw[color=black, very thick](15,0) circle (1); 
\draw[color=black, very thick](5,-5) circle (1); 

\draw[color=red, very thick, double distance = 12pt] (1,0) -- (4,0)  ;
\draw[color=red, very thick, double distance = 12pt] (6,0) -- (9,0)  ;
\draw[color=red, very thick, double distance = 12pt] (11,0) -- (14,0) ;

\draw[color=red, very thick, double distance = 12pt] (5,-1) -- (5,-4) ;
\draw[color=red, very thick, double distance = 12pt] (5,1) -- (5,2) ;
\draw[color=red, very thick, double distance = 12pt] (10,1) -- (10,2) ;
\draw[color=blue, double distance = 15pt] (14.7,1) arc (20: 160: 7.6); 
\draw[color=blue, double distance = 15pt] (6,-5) arc (-90: -39: 11); 

\draw (0,0) node{$1$} ;
\draw (5,0) node{$2$} ;
\draw (10,0) node{$3$} ;
\draw (15,0) node{$4$} ;
\draw (5,-5) node{$5$} ;
\draw (7.5,7.7) node{$\text{loop 1}$} ;
\draw (12,-5) node[below]{$\text{loop 2}$} ;
\draw (5,2) node[above]{$\text{cilium 1}$} ;
\draw (10,2) node[above]{$\text{cilium 2}$} ;
\end{tikzpicture}
\caption{A LVR graph In the case $p=2$  with five vertices colored in black, four propagators colored in red, two loops colored in blue, two cilia and one broken face colored in gray.} 
\label{Fig3a}
\end{figure}

The amplitude of a LVR graph with $k$ cilia expands in trace invariants as:
\bee \label{amplitudetraceinvariant}
 {\cal A}_{(G,T)}^k(\lambda,N,J) =\sum_{\pi\in\Pi_{k}}{A}_{(G,T)}^{\pi}(\lambda,N) \; 
 \Tr_{\pi}\{J,J^{\dagger}\} \; ,
\ee
with $ Tr_{\pi}\{ J,J^{\dagger}\}=\Tr [J J^{\dagger}]^{k_{1}} \cdots \Tr [JJ^{\dagger}]^{{k}_{p}}$ is a trace invariant defined by an 
increasing sequence of $|\pi|$ integers $0<{k}_{1}\leq \cdots\leq {k}_{|\pi|}$ such that ${k}_{1}+ \cdots + {k}_{|\pi|}={k}$,
and the definition of ${A}_{(G,T)}^{\pi}(\lambda,N)$ is
\begin{definition}[Weingarten definition]\label{prop:ampliWeingarten}
\bea
&&{A}_{(G,T)}^{\pi}(\lambda,N)=
\frac{(-\lambda)^{e(G)}N^{v(G)-e(G)}}{v(G)!}
\int dw_{T} \partial_{T} \int d\mu_{C\{x_{ij}^\cT\}}(M)\\&&\nonumber
\sum_{\substack{\tau,\sigma\in{\mathfrak S}_{k} \\ C(\sigma)=\pi}}\;
\sum_{1\leq p_{1},\dots,p_{k}\leq N}\text{Wg}(\tau\sigma^{-1},N)
\prod_{1\leq m\leq f(G)-b(G)}\Tr\Big[{\mathfrak X}^{m}\Big]\prod_{1\leq l\leq k}
{\mathfrak X}^{l}_{p_{\tau(l) }p_{\zeta(l)}}\label{amplitudeWg} \; .
\eea
where ${\mathfrak X}$ is defined by \eqref{Xnew} and \eqref{Ynew}.
\end{definition}
 
If the LVR graph $(G,T)$ is reduced to a tree we use the shorthand notation ${A}_{T}^{\pi}(\lambda,N) $ instead of 
${A}_{(T,T)}^{\pi}(\lambda,N)$. 

\begin{proposition}\label{structure:prop}
The \emph{scalar cumulants} of order $2{\mathcal K}$  can be written as a sum over  partitions of ${\mathcal K}$ and over two permutations of ${\mathcal K}$ elements:
\begin{equation}\label{structure:eq} 
{\mathfrak{K}}^{\mathcal K} (\lambda,N) \!=\!\sum_{\pi\in\Pi_{{\mathcal K}} }{\mathfrak{K}}^{\mathcal K}_{\pi}(\lambda,N)
\sum_{\rho,\sigma\in\mathfrak{S}_{{\mathcal K}}}
\prod_{1\leq l\leq {\mathcal K}}\delta_{d_{l},a_{\rho\tau_{\pi}\sigma^{-1}(l)}}\delta_{c_{l},b_{\rho\xi_{\pi}\sigma^{-1}(l)}} ,
\end{equation}
where ${\mathfrak{K}}^{\mathcal K} (\lambda,N)$ is defined by Definition \ref{cumulantsdef} and
where $\tau_{\pi}$ and $\xi_{\pi}$ are arbitrary permutations such that $\tau_{\pi}(\xi_{\pi})^{-1}$ has a cycle structure corresponding to the partition $\pi$. Then the \emph{scalar cumulants} $ {\mathfrak{K}}^{\mathcal K}_{\pi}(\lambda,N) $ are given by the expansion:
\begin{equation}
{\mathfrak{K}}^{\mathcal K}_{\pi}(\lambda,N)=\sum_{T \text{ LVR tree with $k$ cilia}}{\cal A}^{\pi}_{T}(\lambda,N)\label{Kpitree} \; .
\end{equation}
\end{proposition}
\prf This proposition was proved in \cite{GuKra}. \qed
 
\section{Results}
\label{section4}
We recall that all the results of this section  {\bf are valid only for $1\le {\mathcal K}\le {\mathcal K}_{\max}$, where $ {\mathcal K}_{\max}$ is fixed}. 

As our main Theorem we state an analyticity result and a Borel summability for the constructive expansion of cumulants.
\begin{theorem}
\label{th2}
Let $1\le {\mathcal K}\le {\mathcal K}_{\max}$. The expansion 
\begin{equation}
{\mathfrak{K}}^{\mathcal K}_{\pi}(\lambda, N)=\sum_{T\text{ LVR tree with $2{\mathcal K}$ cilia }}{\cal A}^{\pi}_{T}(\lambda, N)\;
\end{equation}
defines an analytic function of $\lambda\in{\cal C}$. Moreover, each term in this sum is bounded as:
\begin{equation}\label{treecumulantsbound}
\big|{\cal A}^{\pi}_{T}(\lambda, N)\big|\leq\frac{N^{2-|\pi|}|\lambda|^{e(T)}\,({\mathcal K}!)^2 \, 2^{2{\mathcal K}}}
{(\cos\frac{\arg\lambda}{(p-1)})^{2e(T)+{\mathcal K}}\,v(T)!} \; ,
\end{equation}
where $|\pi|$ is the number of integers in the partition $\pi$. This expansion reads:
\begin{equation}
{\mathfrak{K}}^{\mathcal K}_{\pi}(\lambda, N)=\sum_{\substack{G\text{ labeled ribbon graph with ${\mathcal K}$ cilia,}\\\text{ broken faces corresponding to $\pi$, and $e(G)\leq n$}}}
\frac{(-\lambda)^{e(G)}N^{\chi(G)}}{v(G)!}
+{\cal R}^{\mathcal K}_{\pi,n}(\lambda,N) \; ,
\end{equation}
where ${\cal R}^{\mathcal K}_{\pi,n}(\lambda,N)$ is a sum over LVR graphs with  ${\mathcal K}$ cilia, at least $n+1$ edges and at most $n+1$ loop edges.
This remainder  is uniformly analytic for $\lambda \in {\cal C}$, and
it obeys the bound, for $\sigma$ constant large enough,
\bea
\Big|{\cal R}^{\mathcal K}_{\pi,n}(\lambda,N)\Big|&\leq& \; \sigma^n \, N^{2-|\pi|}\: [(p-1)n]! \; |\lambda|^{n+1} .
\eea
So it obeys the theorem stated in the Appendix of this article  (Borel-LeRoy-Nevanlinna-Sokal) with $q \to p-1$, $z \to \lambda$,
$\omega\to N$ whenever $N\in{\mathbb N}^*$.
\end{theorem}

In \cite{GuKra}, which is concerned by the LVE, i.e. $p=2$, the terms in the action  
is given by
$N \Tr \bigg[J\big( 1 -i  \sqrt{\frac{\lambda}{N}}A\big)J^{\dagger}\bigg]$.
Differentiation with respect to the field $A$, arising from the Gaussian integration, leads to the creation of a ciliated vertex, and the source terms appear as a pair $J J^{\dagger}$.

Now we are concern with the LVR, i.e. $p \leq 2$, this requires more work.

\begin{figure}[!ht]
\begin{center}
{\includegraphics[width=9cm]{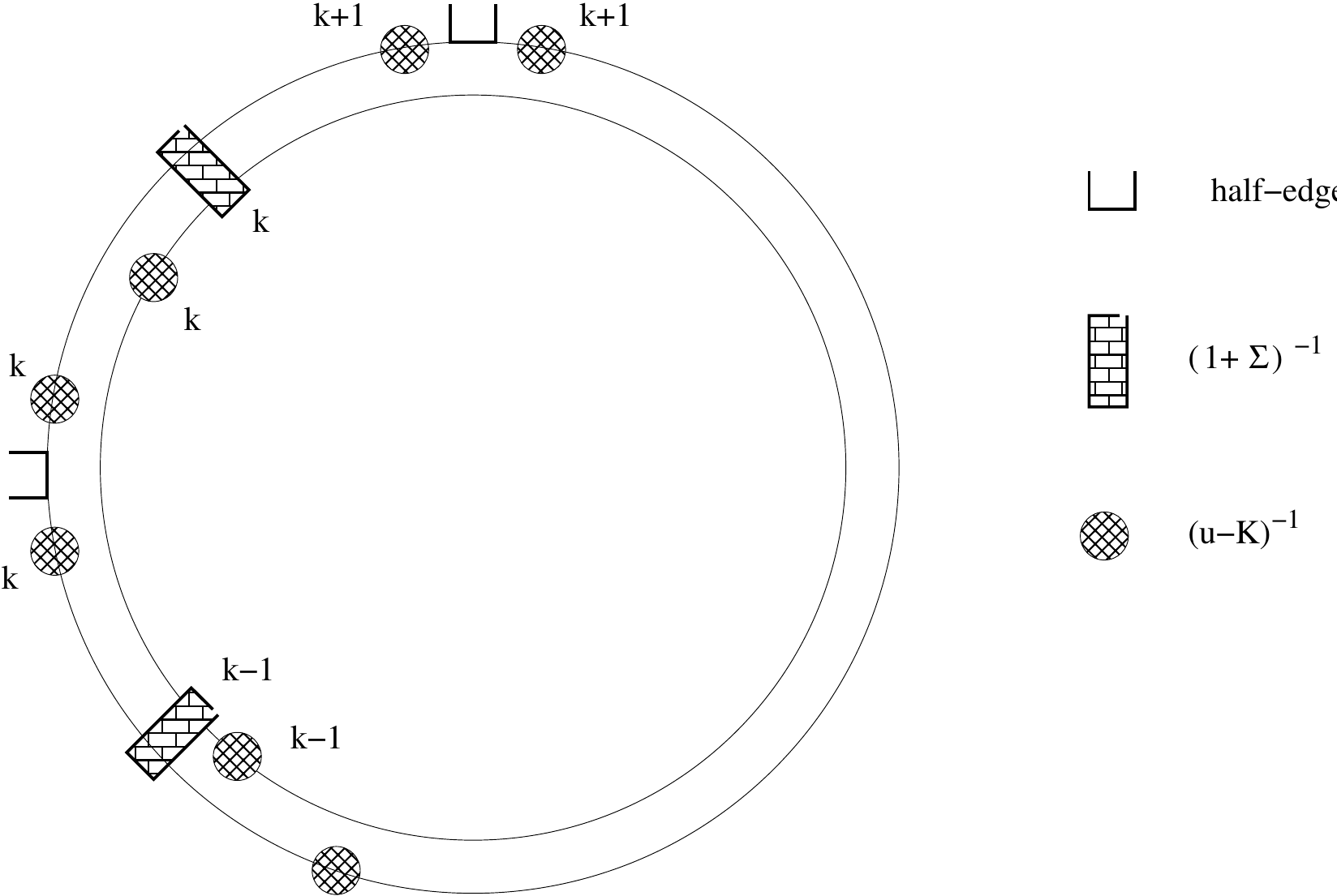}}
\end{center}
\caption{A vertex with some of its corner operators (courtesy from \cite{KRS1}). The label $k$ indicates the corresponding contour variable. The upper left corner between the two half-edges $\sqcup$ symbols contains three $(u-K)^{-1}$ operators with indices $k$, $k$ and $k+1$.}
\label{corneropfig}
\end{figure}

\begin{definition}
\label{defS}
\bea
\Sigma (\lambda, X) &:=& \sum_{k=0}^{p-1} A^k(X)\otimes A^{p-1-k}(X) \\
{\cal S}(J,J^\dagger, N)&:=&N \Tr \bigg[J \bigg(\bbone_\otimes   + \Sigma (\lambda, X) \bigg) J^{\dagger}\bigg]
\eea
\end{definition}

\begin{proposition}
\label{proS}
A half-edge $\sqcup$ can be marked $M$, $M^\dagger$, $J$ or $J^\dagger$.
\end{proposition}

\prf
The proof is given in Section \ref{informal}. 

\qed

\begin{lemma}
Let  $1\le {\mathcal K}\le {\mathcal K}_{\max}$.
In all expressions involving the LVR trees,
we have 
\bea
\bigg( 1 - i   \sqrt{\frac{\lambda}{N}} A \bigg)^{-1}  (J J^{\dagger})^{\eta_c}&=&
\bigg(\bbone_\otimes   + \Sigma(\lambda, M) \bigg)^{-1} (J J^{\dagger})^{\eta_c}.
\eea 
\end{lemma}

\prf 
Taking into consideration Lemma 2.1 and equation (2-30) from the algebraic rule (2.20) and the functional equation (2.12) from \cite{KRS}, which states:
\bea \partial_{X} A(X)&=& \hskip-.2cm\bigg(\bbone_\otimes   + \Sigma(\lambda, M) \bigg)^{-1}=\bigg(\bbone_\otimes   +\lambda \Sigma(\lambda, X) \bigg)^{-1}\hskip-.2cm \text{, where } X=MM^{\dagger} , \nonumber\\
A(X)&=& XT_{p} \bigg(-\frac{\lambda}{N^{p-1}}  X  \bigg),
\eea
the sources terms emerge from the differentiation of ${\cal S}(J,J^\dagger, N)$ with respect to $M$ or $M^{\dagger}$.
We have
\bea \partial_{M} {\cal S}(J,J^\dagger ,N) &=&N \big(\bbone_\otimes  + \Sigma(\lambda, M) \big)^{-1} J^{\dagger}, \\
\partial^{\dagger}_{M} {\cal S}(J,J^\dagger, N) &=&NJ + N \big(\bbone_\otimes   + \Sigma(\lambda, M) \big)^{-1} 
M (M^{\dagger})^{-1} J^{\dagger}\nonumber\\
&-& MM^{\dagger} T_{p} \bigg(-\frac{\lambda}{N^{p-1}} MM^{\dagger}  \bigg)M (M^{\dagger})^{-2} J^{\dagger}.
\label{S25}
\eea
Now we validate such a proposition {\it for  $1\le {\mathcal K}\le {\mathcal K}_{\max}$ by recurrence}.

\begin{itemize}

\item{$ {\mathcal K}=1$}

From \cite{KRS} and \cite{KRS1} we know that the broken face is equal to $1$ and in this broken face there are 
two half-edges $\sqcup$, one marked $J$ and one marked $J^\dagger$ joined by a single line of sources $JJ^\dagger$ comparable to that of Figure \ref{corneropfig}. 

Now $\partial_{M} {\cal S}(J,J^\dagger, N)$ is proportional to $J^{\dagger}$ with a coefficient that does not depend on $J$. 
$\partial_{M}^{\dagger} {\cal S}(J,J^\dagger, N)$ is a sum of $NJ$ and a sum proportional to $J^{\dagger}$ with a coefficient that does not depend on $J$ and $J^{\dagger}$.

Therefore all the other half-edges $\sqcup$ in the broken face are marked $M$ or $M^\dagger$.
 From that we can conclude.

\item{$ [{\mathcal K} \leq s] \Rightarrow   [{\mathcal K}=s+1]$}

Suppose ${\mathcal K} \leq s$ is true. The case ${\mathcal K}=s+1$ can be decomposed into 

\subitem -  one broken face with ${s}_{1}\leq  s$ cilia and the rest of faces (whether they are broken or unbroken).
This case is easy since the sources $J$, $J^\dagger$ are trapped by their broken faces (see Figure  \ref{Fig3a}).

\subitem - or one broken face with ${s}_{1} = s+ 1$ cilia  and the rest of unbroken faces. 
This case is the hard part of the Lemma. There are $s+1$ $J$ who must come from  $s_1\geq s+1$ derivates 
$\partial_{M}^{\dagger}$ since $\partial_{M}$ have no $J$. But there are  $s+1$ $J^\dagger$ (those which come
from $\partial_{M}$ and $\partial_{M}^{\dagger}$).  Now the number of the $J$ and number of the $J^\dagger$ 
must be equal  no matter where the gluing of the half-edges $\sqcup$ are, since the total number the sources $J$, $J^\dagger$ are trapped by the single broken face, 
the number of $J$ are equal to $s+1$ and also the number of $J^\dagger$ must be $s+1$. From this it is easy to conclude.

\end{itemize}
\qed

\section{Summary of \texorpdfstring{\cite{KRS}}{KRS} and \texorpdfstring{\cite{KRS1}}{KRS1}}
\label{secder}

We need now to compute $ \partial_\cT^M \cR_n$. This will be relatively easy since $\cR_n$ is a product of resolvents
of the $\frac{1}{u-X}$ type. 
Since trees have arbitrary coordination numbers we need a formula
for the action on a vertex factor $\cR^i $ of a certain number $r^i=q^i+ \bar q^i$ of
derivatives, $q^i$ of them of the $\frac{\partial }{\partial M_i }$ type and $\bar q_i $ of the 
$\frac{\partial }{\partial M_i^\dagger} $ type.

Let us fix a given loop vertex and forget for a moment the index $i$. We need
to develop a formula for the action of a differentiation operator 
$\frac{\partial^r }{\partial M_1 \cdots \partial M_q  \partial M^\dagger_1  \cdots \partial M^\dagger_{\bar q} } $ on $\cR = \Tr \frac{1}{v_1-X}  \otimes  \Tr \frac{1}{v_2-X}  $.

To perform this computation, we first want to know on which of the two traces (also simply called ``loops") of a loop vertex the differentiations act.
Therefore, we add to any oriented tree $\cT$ of order $n$ a collection of $2(n-1)$ indices $s_e$. Each such index takes value in $\{1,2\}$, 
and specifies at each end $e$ of an edge of the tree whether the field derivative for this end hits the $ \Tr \frac{1}{v_1-X}$ loop or the  
$ \Tr \frac{1}{v_2-X}$ loop. 
There are therefore exactly $2^{2(n-1)}$ such decorated 
oriented trees for any oriented tree. Unless otherwise specified
in the rest of the paper we simply use the word ``tree" for an oriented decorated tree with these additional $\{s \}$ data.
Similarly the set $\cT_n$ from now on means the set of oriented \emph{decorated} trees at order $n$. 

Knowing the decorated tree $\cT$, at each vertex we know how to decompose the number of differentiations
acting on it according to a sum over the two loops of the number of differentiations on that loop, as $q = q_1 + q_2$, $\bar q= \bar q_1 + \bar q_2$.
Hence we have the simpler problem to compute the differentiation operator  
$\frac{\partial^r }{\partial M_1 \cdots \partial M_{q}  \partial M^\dagger_1  \cdots \partial M^\dagger_{\bar q}} $ on 
a \emph{single} loop $\Tr \frac{1}{v-X} $.

We shall use the symbol $\sqcup$ to indicate the place where the indices of the derivatives act\footnote{The symbol $\sqcup$
instead of $\otimes$ will hopefully convey the fact that these derivatives are half propagators for the LVR. The edges of the LVR always glue two $\sqcup$ symbols together.}.
For instance we shall write 
\bee  \frac{\partial }{\partial X}  \frac{1}{v-X} =    \frac{1}{v-X} \sqcup \frac{1}{v-X} .
\ee
To warm up let us compute explicitly some derivatives (writing $\partial_M$ for  $\frac{\partial }{\partial M}$) \:
\begin{eqnarray}
\nonumber
\partial_M \Tr \frac{1}{v-X}  &= &\Big[\Tr\frac{1}{v-X} \sqcup M^\dagger  
 \frac{1}{v-X} \Big] 
\label{eqder1}\\
\partial_{M^\dagger}  \Tr \frac{1}{v-X} &= &\Big[\Tr\frac{1}{v-X} M \sqcup \frac{1}{v-X}  \Big]  .
\label{eqder2}
\end{eqnarray}
Induction is clear: $r=q+\bar q$ derivatives create insertions of $\sqcup M^\dagger  $ and of 
$M \sqcup$ factors in all possible cyclically distinct  ways but they can 
also create double insertions noted $ \sqcup \sqcup $ when a $M^\dagger  $
or $M$ numerator is hit by a derivative. For instance, at second order 
we have:
\begin{eqnarray}
\nonumber
\partial_M^{\textcolor{white}\dagger}  \partial_{M^\dagger}  \Tr \frac{1}{v-X}  
&= &\Tr\Big[\frac{1}{v-X} M \sqcup \frac{1}{v-X} \sqcup M^\dagger  
 \frac{1}{v-X} \Big] 
 \nonumber\\
 &+ &\Tr\Big[\frac{1}{v-X} \sqcup M^\dagger  
 \frac{1}{v-X} M \sqcup \frac{1}{v-X} \Big]  
 \nonumber\\
 &+ &\Tr\Big[\frac{1}{v-X} \sqcup \bbone \sqcup
 \frac{1}{v-X} \Big]  .
 \label{eqder3}
\end{eqnarray}
Remark the last term in which the second derivative hits the numerator created by the first. 
Since $X =M M^\dagger$ the outcome for a $q$-th order partial derivative, 
is a bit difficult to write, but the combinatorics is quite inessential for our future analyticity bounds.
The Fa\`a di Bruno formula allows to write this outcome as as sum over a set  $\Pi^{q, \bar q}_r$ of Fa\`a di Bruno  terms 
each with prefactor 1:
\bee \frac{\partial^r }{\partial M_1 \cdots \partial M_q  \partial M^\dagger_1  \cdots \partial M^\dagger_{\bar q} }\Tr  \frac{1}{v-X} 
= \sum_{\pi \in \Pi^{q, \bar q}_r} \;\Tr\Big[O^\pi_0 \sqcup O^\pi_1 \sqcup \cdots \sqcup O^\pi_r \Big]. \label{faasum}
\ee
In the sum \eqref{faasum} there are exactly $r$  symbols $\sqcup$, separating $r+1$ corner operators $O^\pi_c$.
These corner operators can be of four different types, either resolvents $\frac{1}{v-X}$, $M$-resolvents $\frac{1}{v-X} M $, $M^\dagger$-resolvents 
$M^\dagger \frac{1}{v-X}  $,
or the identity operator $\bbone$. We call $r_\pi$, $r^M_\pi$, $r^{M^\dagger}_\pi$ and $i_\pi$ the number of corresponding operators 
in $\pi$. We shall need only the following facts.
\begin{lemma} We have
\bee \vert \Pi^{q, \bar q}_r \vert \le   2^r  r!, \quad
r_\pi   =   1  +  i_\pi    ,  \quad   r^M_\pi  +r^{M^\dagger}_\pi = r - 2 i_\pi .
\ee \label{faacomb}
\end{lemma}
\proof The lemma is proved in \cite{KRS} under the name of Lemma 3.1.  \qed

Applying \eqref{faasum} at each of the two loops of each loop vertex, we get for any decorated tree $\cT$
\bee  \partial_\cT^M   \cR_n = \prod_{i=1}^n \Big\{ \prod_{j= 1}^2 \Big[ \sum_{\pi_j^i \in \Pi_{r^i_j}^{q^i_j, \bar q^i_j  }  } 
\Tr  \big( O^{\pi_j^i}_0 \sqcup O^{\pi_j^i}_1 \sqcup \cdots \sqcup O^{\pi_j^i}_{r^i_j} \big) \Big]  \Big\}\label{faafullsum}
\ee
where the indices of the previous \eqref{faasum} are simply all decomposed into indices for each loop $j=1, 2$
of each loop vertex $i = 1 ,\cdots , n $.

We now need to understand the adhesion of the $\sqcup$ symbols. 
Knowing the decoration of the tree, that is, the $2(n-1)$ indices $s_e$, we know exactly for which edge of the decorated
tree it connects to which loops. In other words, the \emph{decorated}
tree $\cT_n$ defines a particular \emph{forest} on the $2n$ loops of the $n$ loop vertices. This forest having $n-1$ 
edges must therefore have exactly $n+1$ connected components, each of which is a tree but on the $2n$ loops. We call these trees the \emph{cycles} $\cC$ of the tree, since as trees they have a single face.

\begin{figure}[!ht]
\begin{center}
{\includegraphics[width=12cm]{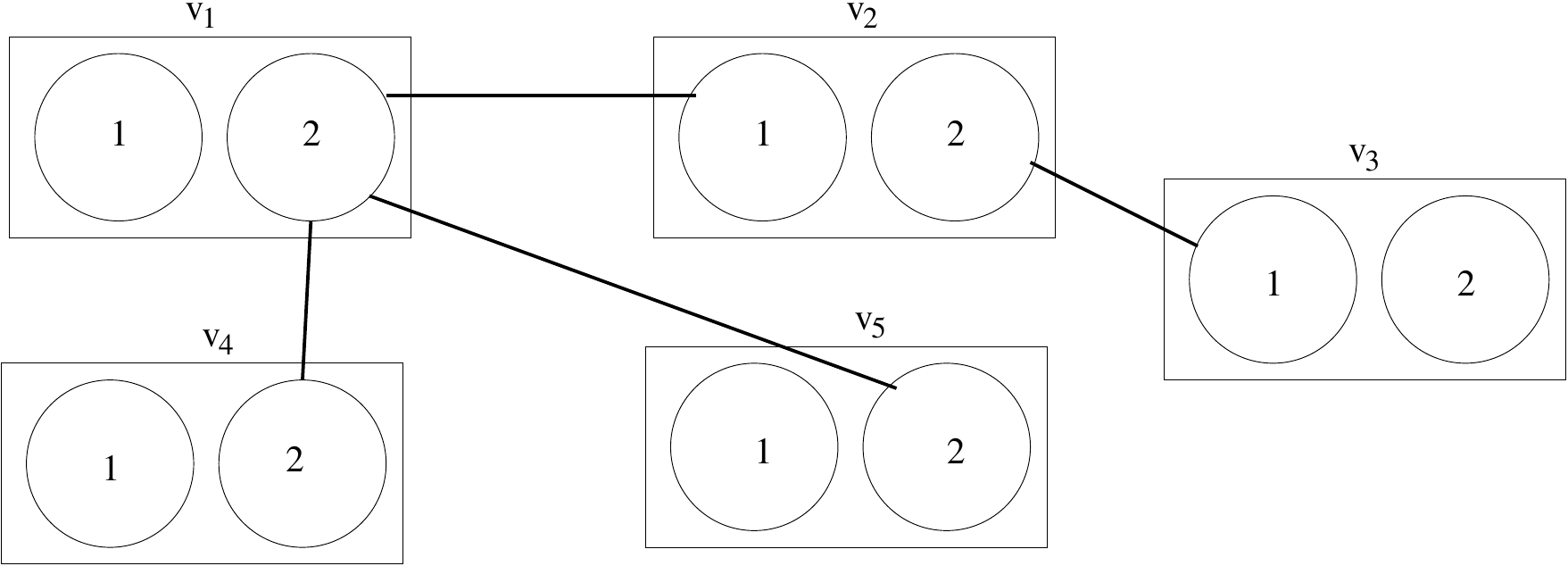}}
\end{center}
\caption{A tree of $n-1$ lines on $n$ loop vertices (depicted as rectangular boxes, hence here $n=5$) 
defines a forest of $n+1$ connected components or cycles $\cC$ on the $2n$ elementary loops, since
each vertex contains exactly two loops. To each such cycle
corresponds a  trace of a given product of operators in the LVR.}
\label{forest}
\end{figure}

Now a moment of attention reveals that, if we fix a particular choice $\{ \pi_j^i \}$ in the expansion of \eqref{faafullsum} obtained by the action of $ \partial_\cT^M $ on $\cR_n$, 
the $\sqcup$ symbols are summed with indices forced to coincide along the edges of the tree, and
we simply glue the $2n$ traces of \eqref{faafullsum} into $n+1$ traces, one for each cycle $\cC$ of the decorated tree $\cT$.
This is the fundamental characteristic of the LVR \cite{rivasseau2007constructive}.
Each trace acts on the product of all corner operators $O_c$ cyclically ordered as obtained by turning the cycle $\cC$. Hence 
we obtain, with hopefully transparent notations,
\bee \partial_\cT^M   \cR_n = \prod_{i=1}^n \Big\{ \prod_{j= 1}^2 \Big[  \sum_{\pi_j^i \in \Pi_{r^i_j}^{q^i_j, \bar q^i_j  }  }  \Big] \Big\} \prod_{\cC}
 \Bigl[ \Tr  \prod_{c  \; \circlearrowleft \; \cC} O_c \Big] .  \label{faafullsumcycl}
\ee
We now bound the associated tree amplitudes of the LVR.

\begin{lemma} \label{mainamp}
For any $\epsilon>0$, there exists $\eta_{\epsilon}>0$ and a constant  $K>0$  such that for any 
tree $\cT$ with $n$ vertices the amplitude $ A_{\cT}(\lambda,N)$ is analytic 
in $\lambda$ in the pacman domain $P(\epsilon, \eta_\epsilon)$
and satisfies in that domain to the \emph{uniform} bound in $N$
\begin{align} \label{mainbou}
 | A_\cT (\lambda,N )| \le K^n \vert \lambda \vert ^{\kappa_p n} \prod_{i=1}^n r_i ! 
\end{align}
where $r_i \ge 1$ is the coordination of the tree $\cT$ at vertex $i$. 
\end{lemma}
\proof This Lemma is proved in \cite{KRS}.
\qed

Then a loop vertex of the theory can be pictured as in Figure \ref{matrixlve2},
where the cilium and each derived leaf bear a factor $\sqcup$, each edge bears a (tensor) resolvent $R$
and each ordinary leaf bears a factor $A_p$.

In \cite{KRS} the following equation holds
\bee 
Z(\lambda, N) = \sum_{n = 0}^\infty \frac{1}{n!}\,\sum_{\cF\in \gF_n} N^{-\vert \cF \vert}  \int dw_\cF 
\int \{dt du dv\} \Phi_n  \int d\mu_{C(x)} \partial_\cF^M   \cR_n
 \Big|_{x_{ij} = x_{ij}^\cF (w)} .
  \label{LVE4}
\ee

Similarly a LVR tree is obtained by gluing $n$ such loop vertices through along $n-1$ pairs of 
glued $\sqcup$ factors, see Figure \ref{matrixlve1}.
Beyond the tree, additional cycles between the loop vertices can of course exist but they are hidden in the 
functional integral $ \int dw_\cF  \int d\mu_{C(x)} $  in \eqref{LVE4}. 

\section{Proof of Theorem \ref{th2}.}
\label{informal}
\begin{figure}[!ht]
\begin{center}
{\includegraphics[width=7cm]{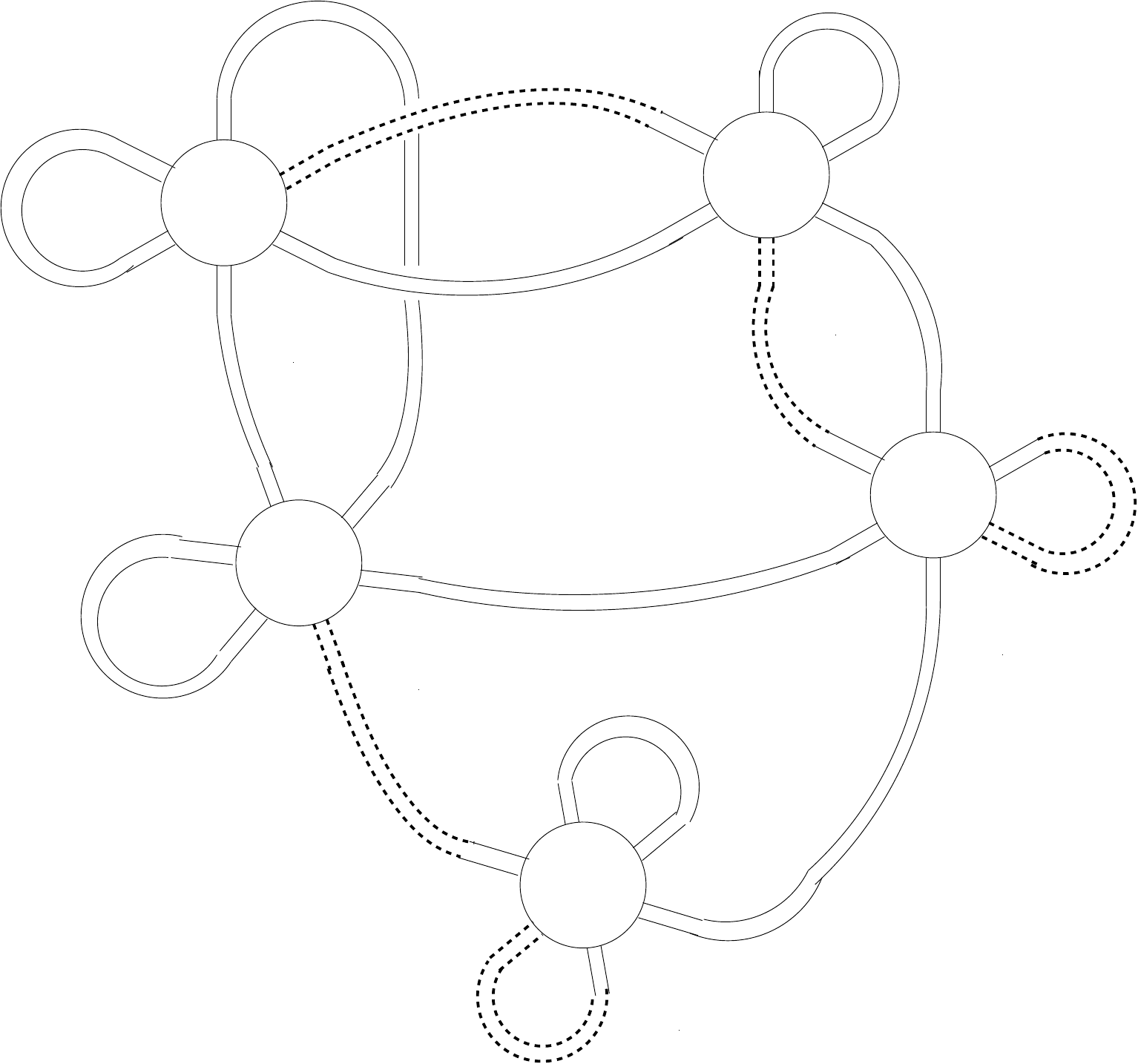}}
\end{center}
\caption{A Feynman graph of the $\Tr (M^\dagger M )^3$ theory. All five vertices are 6-valent. One $M^\dagger$
field per vertex leads to a dotted line (for better visibility we showed as dotted the hooked point plus a large fraction of the propagator). 
We define in this case two connected components, namely two loop vertices.} 
\label{feynman1}
\end{figure}

First, we explain how to visualize  the LVR vertices associated to a Feynman graph with sources $J$, $J^\dagger$.
Taking an ordinary connected Feynman graph, we draw,
at every vertex of the graph, one selected half-edge marked $M^\dagger$ 
 as a dotted half-line. The set of edges which are so dotted then defines a subset of connected components, each of which has a \emph{single loop}. They are the LVR \emph{vertices} associated
to this graph (see Figure \ref{feynman1}). 

\begin{figure}[!ht]
\begin{center}
{\includegraphics[width=7cm]{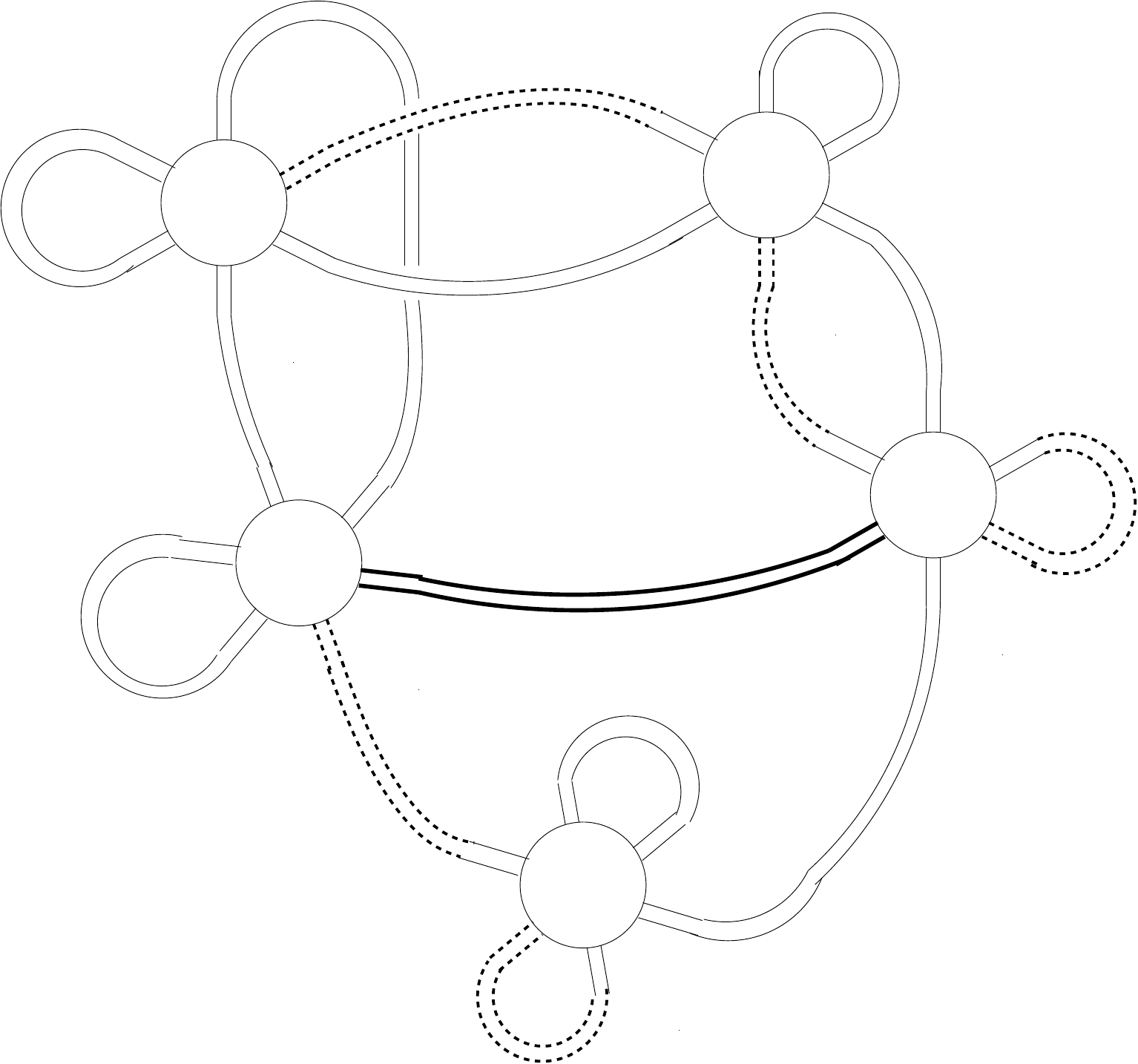}}
\end{center}
\caption{Adding a tree line in boldface between the two loop vertices
gives one of the Feynman graph contributions to the LVR tree made of two loop vertices
joined by a single edge.}
\label{feynman2}
\end{figure}

Selecting a spanning tree between these vertices through the BKAR formula, is like dividing each
Feynman graph built around these LVR vertices into as many pieces as there are 
of \emph{spanning trees} between them. Each piece is then attributed to the corresponding LVR tree
(see Figure \ref{feynman2}).

\begin{figure}[!ht]
\begin{center}
{\includegraphics[width=8cm]{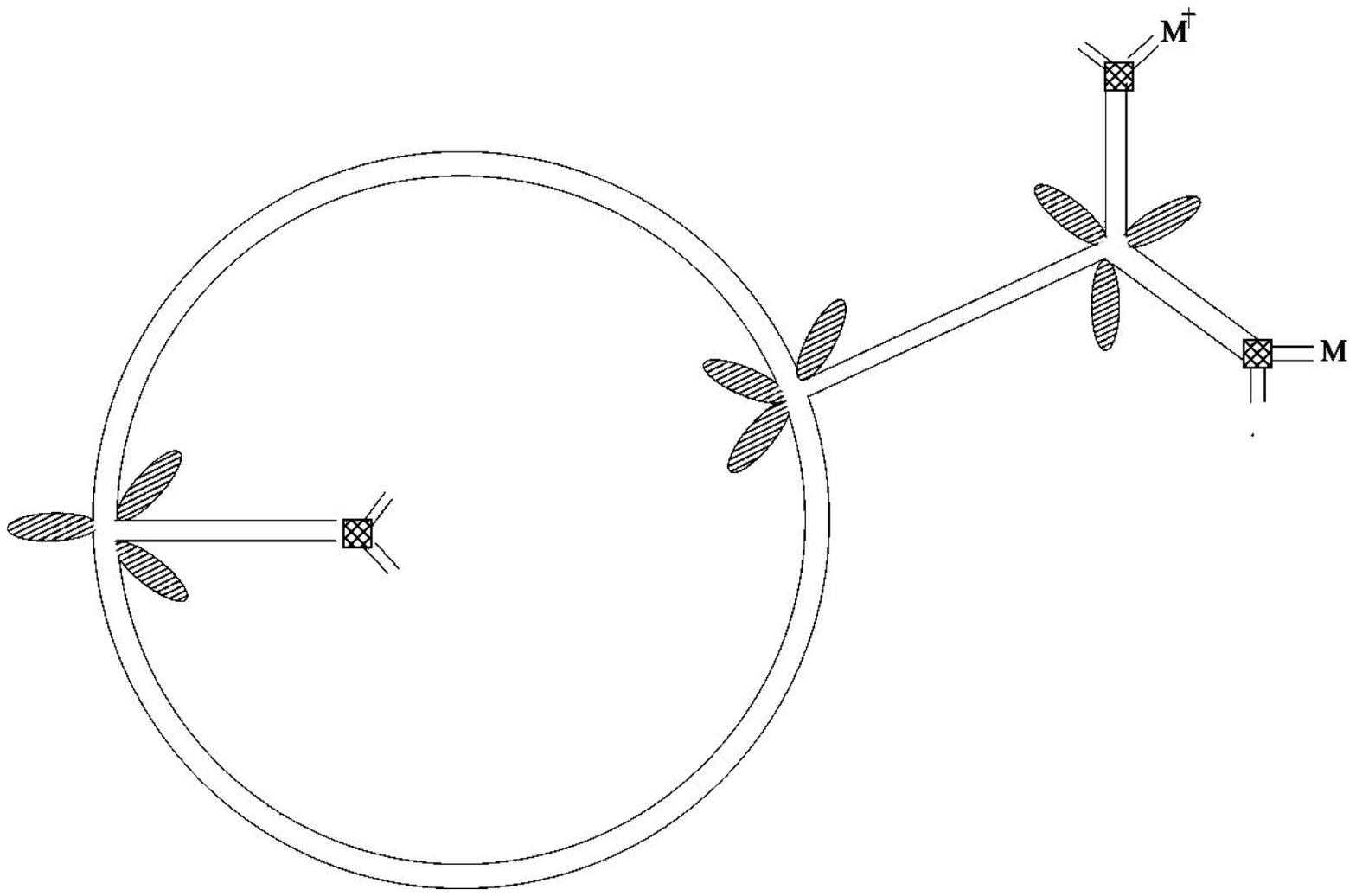}}
\end{center}
\caption{A loop vertex of the theory, bearing 4 derivatives, hence four sources $\sqcup$, where the $\sqcup$
can be marked $M$, $M^\dagger$, $J$ or $J^\dagger$.
We chose $p=5$, hence all vertices are 6-valent. Hatched ellipses represent $A_p$ insertions, ribbon edges represent resolvents
(there are five such resolvents in this graph) and squares represent derived leaves 
which can be of three different types $M \sqcup$, $\sqcup M^\dagger$ or $\sqcup \bbone \sqcup$, where the $\sqcup$
can be marked $M$, $M^\dagger$, $J$ or $J^\dagger$.
In the case pictured, we have three squares because
two derivatives acted on the same $M^\dagger$ factor.}
\label{matrixlve2}
\end{figure}

Conversely if we start from a given LVR tree and want to picture the whole set 
of (pieces of) Feynman graphs that it sums, we introduce
a symbol, such as a hatched ellipse, to picture the sum of all p-ary trees in the 
generating $A_p$ function. A loop vertex of the theory can be then pictured as in Figure \ref{matrixlve2},
where the cilium and each derived leaf bear a factor $\sqcup$, each edge bears a (tensor) resolvent $R$
and each ordinary leaf bears a factor $A_p$.
\begin{figure}[!ht]
\begin{center}
{\includegraphics[width=11cm]{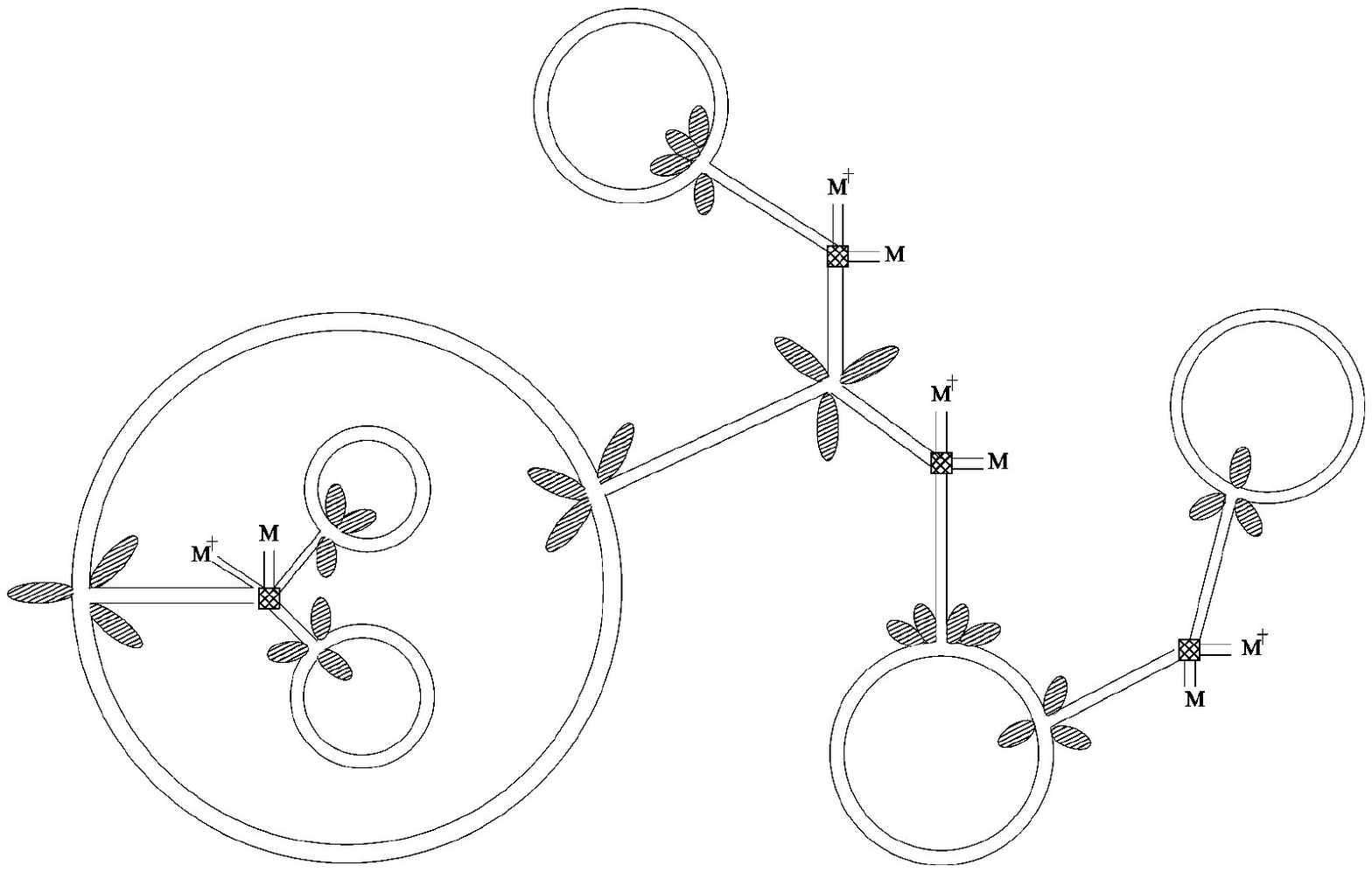}}
\end{center}
\caption{A tree of the loop vertex expansion. It is made of  six loop vertices, 
joined by four edges $MM^\dagger$, each bearing a square, which corresponds to the gluing of two half-edges $\sqcup$,
one marked by $M$ and the other marked by $M^\dagger$.
The attentive reader can find seven traces in the drawing.}
\label{matrixlve1}
\end{figure}

Similarly a LVR tree is obtained by gluing $n$ such loop vertices through along $n-1$ pairs of 
glued $\sqcup$ factors, see Figure \ref{matrixlve1}.
Beyond the tree, additional cycles between the loop vertices can of course exist but they are hidden in the 
functional integral $ \int dw_\cF  \int d\mu_{C(x)} $  in \eqref{LVE4}.

\begin{figure}[!ht]
\begin{center}
{\includegraphics[width=12cm]{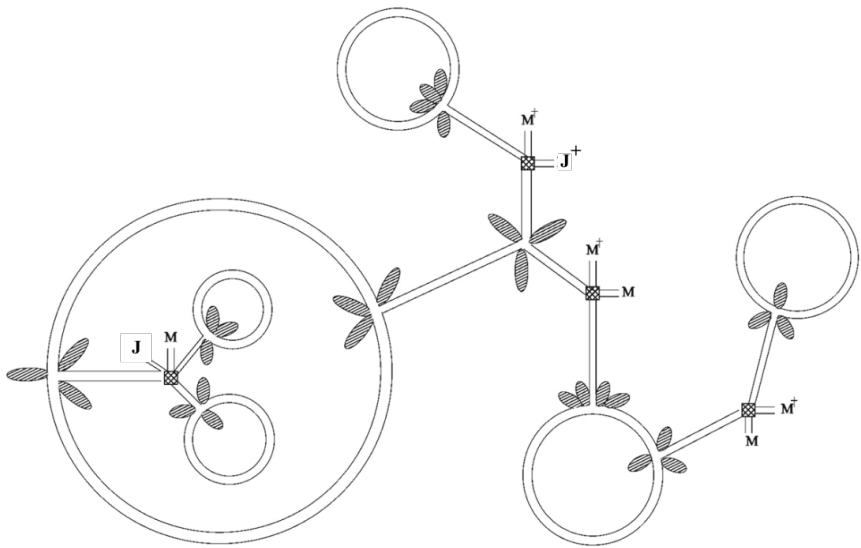}}
\end{center}
\caption{A tree of the loop vertex expansion for ${\mathcal K}=1$, hence (in the language of physicists) for a two point connected Schwinger function. It is made of  six loop vertices, 
joined by four edges each bearing a square, which corresponds to the gluing of two $\sqcup$
of the previous picture, and the attentive reader can find seven traces in the drawing.}
\label{R7n}
\end{figure}

Our first proposition  is a convergent expansion of $\log {\cal Z} (\lambda,N,J)$ as a sum over LVR trees. 
\begin{proposition}
\label{treeexpansion}
For any $\lambda\in {\cal C}$, there exists $\epsilon_{\lambda}>0$ depending on $\lambda$ such that for $\|JJ^{\dagger}\|<\epsilon_{\lambda}$ the 
logarithm of ${\cal Z} (\lambda,N,J) $ is given by the following absolutely convergent expansion:
\begin{equation}
\log {\cal Z} (\lambda,N,J) =   \sum_{T
\,\text{LVR tree}}
{\cal A}_{T}(\lambda,N,J) \;,
\label{treeexpansion:eq}
\end{equation}
\end{proposition}
\proof It is sufficient to combine Theorem 2  in \cite{GuKra} 
with the figure \ref{R7n} of this paper
(which corresponds almost to the figure 9 in
\cite{KRS1}), and (for $p\ge 3$) to add the different loop vertex representation 
including $J$ and $J^\dagger$. 
Remark that the normalization in \cite{GuKra} is different than in \cite{KRS}, 
\eqref{diffnor} takes it into account. 
\qed

\medskip
In order to compare the tree expansion of Proposition \ref{treeexpansion} with the conventional perturbative expansion, it is  necessary to further expand some of the loop edges. The following proposition is obtained by recursively adding loop edges to the LVR trees.

\begin{proposition}
\label{perturbativegenerating:thm}
For any $\lambda \in {\cal C}$, there exists $\epsilon_{\lambda}>0$ depending on $\lambda$ such that for $\|JJ^{\dagger}\|<\epsilon_{\lambda}$:
\bea
\nonumber
&&\hskip-.6cm \log {\cal Z}(\lambda,N,J) = \hskip-.2cm
\sum_{\substack{G\text{ ciliated labeled}\\\text{ ribbon graph}\\e(G)\leq n}}\hskip-.2cm
\frac{(-\lambda)^{e(G)}N^{\chi(G)}}{v(G)}\prod_{f\in b(G)}\Tr (JJ^{\dagger})^{c(f)} 
+  {\cal R}_{n}(\lambda,N,J)\; ,
\label{perturbativeexpansion:eq}
\eea
where $c(f)$ is the number of cilia in the broken face $f$ and the perturbative remainder at order $n$
is a convergent sum over LVR  graphs with at least $n+1$ edges and at most $n+1$ loop edges
\begin{equation}
{\cal R}_{n}(\lambda,N,J)=
\sum_{\substack{(G,T)\text{ LVR graph}\\ e(G)= n+1}}
{\cal A}_{(G,T)}(\lambda,N,J)
+ \sum_{\substack{T\text{ LVR tree}\\ e(T)\geq n+2}}
{\cal A}_{T}(\lambda,N,J) \; .
\end{equation}
\end{proposition}
\prf This proposition was  proved in \cite{GuKra}. 
\qed

Now we cite the main  theorem proved in \cite{Riv1}.
\begin{theorem}[Constructive expansion for the J-cumulant] \label{th1}
Let $1\le {\mathcal K}\le {\mathcal K}_{\max}$, where $ {\mathcal K}_{\max}$ is fixed. There exists 
$\epsilon_{\lambda}>0$ depending on $\lambda$ such that 
$\mathfrak{K}^{{\mathcal K}}(\lambda,N) $ is given by the following absolutely convergent expansion 
\bea
\mathfrak{K}^{{\mathcal K}}(\lambda,N) &=&
{\cal P}^{\mathcal K}_{n}(\lambda,N,J)+{\cal Q}^{\mathcal K}_{n}(\lambda,N,J) + {\cal R}^{\mathcal K}_{n}(\lambda,N,J)\label{treeexpansion00},\\
{\cal P}^{\mathcal K}_{n}(\lambda,N,J)&=&\sum_{\substack{\substack{G\text{labeled ribbon graph }\\ \text{ with ${\mathcal K}$ cilia,}} 
\\ e(G)\leq n}} \hskip-.3cm
\frac{(-\lambda)^{e(G)}N^{\chi(G)}}{v(G)!}\hskip-.1cm \prod_{f\in b(G)} \hskip-.1cm\Tr[(JJ^\dagger)^{c(f)}]\label{treeexpansion01} ,\\
\hskip-.3cm {\cal Q}^{\mathcal K}_{n}(\lambda,N,J)&=&\sum_{\substack{\substack{(G,T)\text{ LVR graph }\\ \text{  with ${\mathcal K}$ cilia}}\\ e(T)= n+1}}{\cal A}^{\mathcal K}_{(G,T)}(\lambda,N,J) ,
 \label{treeexpansion02}
 \\{\cal R}^{\mathcal K}_{n}(\lambda,N,J)&=&\sum_{\substack{\substack{T\text{ LVR tree  }\\ \text{ with ${\mathcal K}$ cilia}}\\ e(T)\geq n+2}}{\cal A}^{\mathcal K}_{T}(\lambda,N,J) .
 \label{treeexpansion03}
\eea
This expansion is analytic for any $\lambda\in {\cal C}$ 
and the remainder at order $n$ obeys, for $\sigma$ constant large enough, the analog of \eqref{BorLeRSum0}  
\bee
| {\cal R}^{\mathcal K}_{n}(\lambda,N,J)| =  \big| \mathfrak{K}^{{\mathcal K}}(\lambda,N,J)-\sum_{m=0}^{n}a_{m}(N,J)\lambda^{m} \big|\le \sigma^n\,
[(p-1)n] ! \, | \lambda |^{n+1} ,\label{BorLeRSum01}
\ee
uniformly in $N\in{\mathbb N}^*$, J such that $\|J^{\dagger}J\|<\epsilon_{\lambda}$.
Therefore  it obeys the theorem stated in the Appendix  of this article (Borel-LeRoy-Nevanlinna-Sokal) with $q \to p-1$, $z \to \lambda$, 
$\omega\to \big\{ N,J \big\}$, whenever $N\in{\mathbb N}^*, \|J^{\dagger}J\|<\epsilon_{\lambda}$.
\end{theorem}
\proof 
This theorem is proved in \cite{Riv1}.
\qed

Then Theorem  \ref{th2}  would follow easily from the proof of Theorem \ref{th1} and \cite{GuKra}.
\qed

\section{Appendix: Borel-LeRoy-Nevanlinna-Sokal theorem}

We recall the following theorem \cite{Sokal,CaGrMa}. 
\begin{figure}[htb]
\[
\begin{array}{cc}
\includegraphics[width=4cm]{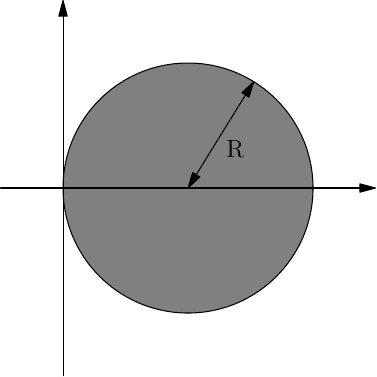}&
\includegraphics[width=6.7cm]{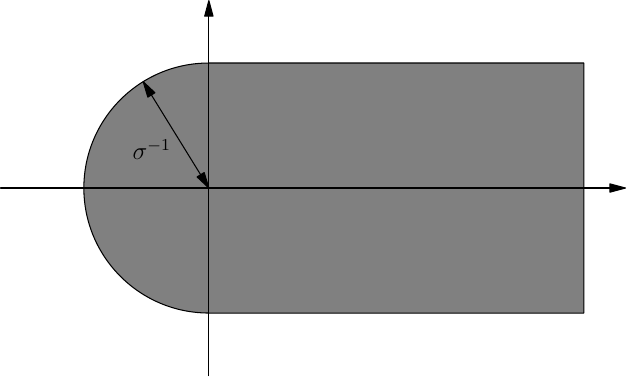}\\
{\cal D}_{R}&\Sigma_{\sigma}
\end{array}
\]
\caption{Domain of analyticity of $F$ and of its Borel transform for $q=1$.}
\label{Borel:pic}
\end{figure}

\begin{theorem}
\label{BorLeRSum}
Let $q\in {\mathbb N^*}$. Let  $F_{\omega}(z)$ be a family of analytic functions on the domain 
\bee
D_R = \{z:\Re z^{-\frac{1}{q}} >(2R)^{-1}\} = \{z:  | z | < (2R)^q \cos^q (\frac{\arg z}{q}) \} 
\label{BorLeRSumbis0}
\ee
depending on some parameter $\omega\in\Omega$, and such that, for some $\sigma \in {\mathbb R}_+$,
\bee
|R_n(z)| =  \big| F_{\omega}(z)-\sum_{m=0}^{n}a_{m}(\omega)z^{m} \big|\le 
\sigma^n
(qn) ! \, | z |^{n+1} \label{BorLeRSum0}
\ee
uniformly in $D_R$ and $\omega\in\Omega$. Then the formal expansion  
\bee\sum_{n=0}^\infty s^{qn} \frac{a_n(\omega)}{(qn)!}\label{BorLeRSumbis1}
\ee
 is convergent for small $s$ and determines a function 
$ B_\omega(s^q )$ analytic in 
\bee
\Sigma_\sigma =\{s :{\text dist}(s ,{\mathbb R}_+ ) < \sigma^{-1} \} 
\label{BorLeRSum1}
\ee 
and such that 
\bee
| B_\omega(s^q )| \le  B \exp\big(\frac{| s |}{R}\big)
\label{BorLeRSum2}
\ee
uniformly for $\Sigma_{\sigma}$ (in \eqref{BorLeRSum2} $B$ is a constant, that is, it is independent of $\omega$). Moreover, setting $t = s^q $,
\bee
 F_{\omega}(z)= \frac{1}{qz}\int_0^\infty B_{\omega}(t)\;
 \Big(\frac{t}{z}\Big)^{\frac{1}{q}-1}  \exp\bigg(-\Big(\frac{t}{z}\Big)^{\frac{1}{q}} \bigg)\; dt
 \label{BorLeRSum3}
\ee
for all $z \in D_R$. Conversely, if $F_{\omega}(z)$ is given by \eqref{BorLeRSum3}, with the above properties for
$B_\omega(s^q)$, then it satisfies remainder estimates of the type \eqref{BorLeRSum0} uniformly, in any $D_r$
such that $0<r <R$, and in  $\omega\in\Omega$.
\end{theorem}
For Theorem \ref{th2}, change $q\to p-1$, $z \to \lambda$, $\omega\to N$ whenever $N\in{\mathbb N}^*$.

\end{document}